\documentclass[12pt]{amsart}
\usepackage{amsfonts, amssymb,amsmath, amsthm, amsxtra, latexsym, amscd}
\usepackage[latin1]{inputenc} % for the swedish symbols äåö etc

\newtheorem{theorem}{Theorem}
\newtheorem{corollary}{Corollary}
\newtheorem{lemma}{Lemma}

\newtheorem{proposition}{Proposition}
\newtheorem{definition}{Definition}
\newtheorem{remark}{Remark}
\newtheorem{example}{Example}

\def\ga{\alpha}

\def\gb{\beta}

\def\gga{\gamma}

\def\gep{\varepsilon}

\def\gv{\varphi}
\def\gp{\Phi}

\def\gl{\lambda}

\def\gO{\Omega}
\def\mf{\mathfrak}
\def\mb{\mathbb}
\def\openA{{{\rm A}\kern-.63em{\rm A}}}
\def\openB{{{\rm I}\kern-.16em {\rm B}}}
\def\openC{{\rm C\kern-.18cm\vrule width.6pt
height 6pt depth-.2pt \kern.18cm}}
\def\openD{{{\rm I}\kern-.16em {\rm D}}}
\def\openE{{{\rm I}\kern-.16em {\rm E}}}
\def\openF{{{\rm I}\kern-.16em {\rm F}}}
\def\openI{{{\rm I}\kern-.2em {\rm I}}}
\def\openH{{{\rm I}\kern-.16em {\rm H}}}
\def\openK{{{\rm I}\kern-.16em {\rm K}}}
\def\openN{{{\rm I}\kern-.16em {\rm N}}}
\def\openP{{{\rm I}\kern-.16em {\rm P}}}
\def\openQ{{\rm Q\kern-.21cm\vrule width.6pt
height 6.2ptdepth-.2pt \kern.21cm}}
\def\openQ{{{\rm Q}\kern-.63em {\rm Q}}}
\def\openR{{{\rm I}\kern-.16em {\rm R}}}
\def\openS{{{\rm S}\kern-.68em{\rm S}}}
\def\openT{{{\rm T}\kern-.30em {\rm T}}}
\def\openZ{{{\rm Z}\kern-.28em{\rm Z}}}

\def\eop{{\hfill\vbox{\hrule height .3pt
      \hbox{\vrule width.3pt height 7pt
      \kern 7pt
      \vrule width .3pt}
      \hrule height .3pt}} \par\bigskip}
\def\meop{\qquad {\vbox{\hrule height .3pt
      \hbox{\vrule width.3pt height 7pt
      \kern 7pt
      \vrule width .3pt}
      \hrule height .3pt}}}

\def\ccd{\cdot\cdot\cdot}
\def\begeqn{\begin{equation}}
\def\endeqn{\end{equation}}
\def\begth{\begin{theorem}}
\def\endth{\end{theorem}}
\def\begprop{\begin{proposition}}
\def\endprop{\end{proposition}}
\def\begcor{\begin{corollary}}
\def\endcor{\end{corollary}}
\def\begdef{\begin{definition}}
\def\enddef{\end{definition}}
\def\beglemm{\begin{lemma}}
\def\endlemm{\end{lemma}}
\def\begexm{\begin{example}}
\def\endexm{\end{example}}
\def\begpr{\begin{proof}}
\def\endpr{\end{proof}}
\def\na{\nabla}
\def\pa{\partial}
\def\bs{\begin{split}}
\def\endsp{\end{split}}
\def\l{\left}
\def\r{\right}

\begin{document}

%%%%%%%%%%%%%%%%%%%%%%%%%%%%%%%%%%%%%

\begin{title}
[$p$-Laplacians and  Hardy-type inequality] {
Degenerate  p-Laplacian operators on H-type groups and
applications to Hardy type inequalities}
\author{Yongyang Jin and Genkai Zhang}
\thanks{Research of Y. Jin supported
by the China State Scholarship and  G. Zhang by the Swedish Research
Council and a STINT Institutional Grant}
\end{title}

\begin{abstract}
Let $\mathbb G$ be a step-two nilpotent group of H-type with Lie
algebra $\mathfrak G=V\oplus \mathfrak t$. We define a class of
vector fields $X=\{X_j\}$ on $\mathbb G$
 depending on
a real parameter $k\ge 1$, and we consider the corresponding
$p$-Laplacian operator $L_{p,k} u= \text{div}_X (|\na_{X} u|^{p-2}
\na_X u)$. For $k=1$ the vector fields  $X=\{X_j\}$ are the left
invariant vector fields corresponding to an orthonormal basis of
$V$, for $k=2$ and $\mathbb G$ being the Heisenberg group they are
introduced by Greiner \cite{Greiner-cjm79}.
 In this paper we obtain the fundamental solution for
the operator $L_{p,k}$ and as an application, we get
 a Hardy type inequality associated with $X$.
\end{abstract}

 \maketitle

\bibliographystyle{plain}

\section{Introduction}
The study of partial differential operators constructed from
non-commutative vector fields satisfying the H\"{o}rmander condition
\cite{Hormander} has had much development. An important class of
such fields, serving as local models,  is that of left-invariant
vector fields on stratified, nilpotent Lie groups with their
associated sub-Laplacians defined by the square-sums of the vector
fields. One of the main tools in the study of the regularity theory
of the sub-Laplacian equation is
 the fundamental
solution; this has been developed in the works of Folland
\cite{Folland1973F.S.} and \cite{Folland1975Arkiv}, Folland and
Stein \cite{Folland-stein1982Hardyspace},  Nagel, Stein and Wainger
\cite{NSW1985vectorfields},  Rothschild  and Stein
\cite{RS1976Acta.Math} and Sanchez-Calle \cite{SC1984Invetion}. In
the papers \cite{CDG1996, HH} the authors studied a class of subelliptic
p-Laplacians on H-type group associated with the left-invariant
vector fields and found the corresponding fundamental solution.

Recently there have been considerable
interests in
studying the sub-Laplacians as square-sums
of vector fields that are not invariant or do not satisfy the
H\"ormander condition. They turn out
to be rather difficult, among the examples of such sub-Laplacians are the
Grushin operators and the sub-Laplacian constructed by Kohn
\cite{Kohn-annmath05}. Those non-invariant sub-Laplacians
also appear naturally in complex analysis.
In the paper  \cite{BGG-JMPA-1998} Greiner, Beals and Gaveau
considered the CR operators $\{Z_j , \bar{Z_j}\}_{j=1}^n$ on $\mb
R^{2n+1}$ as boundary of the complex domain
$$\l\{(z_1 , ... . z_{n+1})\in \mb C^{n+1}:\ \ \text{Im}z_{n+1}> (\sum_{j=1}^n |z_j|^2)^k \r\},$$
where $Z_j =\frac 12({X_j -iY_j}),$  \begeqn X_j =\partial/\partial x_j +2ky_j
|z|^{2k-2}
\partial/\partial t, \ Y_j =\partial/\partial y_j -2kx_j |z|^{2k-2}
\partial/\partial t. \endeqn
and $k$ is a positive integer.  
For $k=1$ these
 vector fields are   left-invariant
on the Heisenberg group $\mathbb R^{2n+1}$
but for $k\ne 1$ there are neither  left nor right-invariant.
The fundamental
solution for their square sum
$$\sum_{j=1}^n Z_j\bar{Z_j} + \bar{Z_j}Z_j$$
is studied  in \cite{BGG-JMPA-1998}. As is well-known,
 the explicit formula of the fundamental solution
is of substantial importance in the study of boundary $\bar\partial$-problem;
see e.g. \cite{Stein-book2}.
In \cite{ZN2003Nonlinear-Analysis}
 Zhang and Niu
studied the
Greiner vector fields on $\mb R^{2n+1}$ for general parameter $k\geq
1$ and found the fundamental solution for the degenerate
$p$-subelliptic operators ${L}_{p,k},$; see Section 2 below. Note
that for non-integral $k$ these vector fields  do not satisfy the
H\"{o}rmander condition and are not smooth.

Heisenberg groups
have  natural generalizations,
namely Carnot groups which are the nilpotent stratified  Lie groups $ \mb G$
with Lie algebras $\mathfrak G=V_1\oplus V_2 \oplus \dots \oplus V_l$
with $[V_i, V_j] \subset V_{i+j}$, with the
sub-Riemannian structure defined by the generating
subspaces $V_1$.  The sub-Laplacian has
generalization to $p$-sub-Laplacian generated
by non-invariant vector fields. The $p$-sub-Laplacian
in this setting
plays important role in the study of quasiregular
maps \cite{HH}. The general theory
in this setup is still not fully developed.

An important subclass of Carnot groups is that of  H-type groups
which were introduced by Kaplan \cite{Kaplan1980FS} as direct
generalizations of Heisenberg groups. In the present paper we will
define a class of vector fields $X$ (see (2.3) below) on H-type
groups generalizing the vector fields (1.1) considered in
\cite{BGG-JMPA-1998} and \cite{ZN2003Nonlinear-Analysis}, and we
find the fundamental solution of the corresponding $p$-Laplacian
with singularity at the identity element. As application we prove a
Hardy type inequality associated to $X$.

Here is a brief review and comparison of our results with those in
the literature. The case of Heisenberg groups with general parameter
$k$
 is done in \cite{ZN2003Nonlinear-Analysis}.
 When   $\mb G$
is a general Carnot group with the sub-Laplacian being invariant
Hardy type inequality has been proved by D'Ambrosio \cite{Ann.Pisa};
see also \cite{DAm-diff.equ} where Hardy type inequalities on
Heisenberg groups are studied. Our  vector fields are however not
invariant and not smooth for non-integral $k$, and our techniques
 are slightly different from theirs. In particular the  computations
in our case  are rather
involved; we use some fine structure of the H-type groups
and we obtain also the best constant for the Hardy type inequality.

The paper is organized as follows. In Section 2 we recall
 some
 basic facts of the H-type group and introduce
the   degenerate p-Laplacian operator $L_{p,k}$
generalizing the invariant sub-Laplacian; Section 3 is
devoted to the proof of the fundamental solution for $L_{p,k};$  In
the final Section 4 we prove the Hardy type inequality associated
with $X$.

\section{H-type groups and a family of vector fields}

We recall that a simply connected nilpotent group $\mb G$ is  of
Heisenberg type, or simply H-type, if its Lie algebra $\mathfrak
G=V\oplus \mathfrak t$ is of step-two, $[V, V]\subset \mf t,$ and if
there is an inner product
 $\langle\cdot\  , \ \cdot \rangle$ in $\mf G$ such that the linear map
 $$J: \mf t \rightarrow  \text{End} (V),$$ defined by the relation
 $$\langle J_{t}(u), v\rangle = \langle t, [u,v] \rangle $$
 satisfies $$J_{t}^{2}=-|t|^{2}\textbf{Id}$$
 for all $t\in \mf t, u,v \in V.$
We denote $m=\text{dim} V$ and $q=\text{dim} \mf t.$

 We identify $\mb G$ with its Lie algebra $\mf G$ via the exponential
  map, $\exp:V \oplus \mf t \rightarrow \mb G.$
The Lie group product is given
by
\begin{equation}
\label{Lie-prod}
(u, t)(v, s)=(u+v, t+s +\frac 12[u, v]).
\end{equation}
Each vector $X\in \mathfrak G$
defines a tangent vector
at any $g$ by differentiating along $g\cdot \exp(tX)$,
namely a left-invariant vector field, denoted also by $X$.
The sub-Laplacian on $\mb G$ is
$$\Delta_{\mb G}=\sum_{j=1}^m X_j^2,$$
where $\{X_j\}$ is an orthonormal basis of $V$.

For
  $g\in \mb G,$ we write $g=(z(g) , t(g))\in V \oplus \mf t,$  and let $K(g)=(|z(g)|^4
+16|t(g)|^2)^{\frac{1}{4}}.$
 In \cite{Kaplan1980FS} Kaplan  proved that there exists a constant $C>0$
 such that the function $$\Phi(g)=C \cdot K(g)^{2-(m+2q)}$$ is a  fundamental
solution for the operator  $\Delta_{\mb G}$ with singularity at the
identity element. We note that $m+2q$ is the homogeneous dimension
of $\mb G$.

In \cite{CDG1996} the authors considered the following
subelliptic $p$-Laplacian
  $$\Delta_p u=\sum_{j=1}^m X_j^* (|\na_{\mb G} u|^{p-2} X_j u)$$ on
H-type group $\mb G$,
where  $\{X_j\}_1^m $ is an orthogonal basis of $V$,   $X_j^*$  is the
formal adjoint of $X_j,$ and $\na_\mb G =(X_1, ... , X_m)$.
For $p=2$ it is the sub-Laplacian above.
They obtained  a remarkable explicit
 formula
 for the fundamental solution
of $\Delta_p$,
 $$ \Gamma_p
=\begin{cases} C_p \ K^{\frac{p-Q}{p-1}}\ , & \quad p\neq Q\\
 C_Q \log \frac{1}{K}\ , & \quad p=Q
\end{cases}$$
As  application, the authors
obtained some regularity results for a class of nonlinear
subelliptic equations.

Motivated by  the work of  Greiner, Beals and Gaveau   \cite{BGG-JMPA-1998},
Zhang and Niu  \cite{ZN2003Nonlinear-Analysis} considered the
following  degenerate p-subelliptic operators on the Heisenberg group
$\mb R^{2n+1}$:
$$L_{p,k} u= \text{div}_{L} (|\na_{L} u|^{p-2} \na_{L} u).
$$
 Here $$ \na_{L}u =(X_1 u ,..., X_n u,Y_1 u ,..., Y_n u), \ \
\text{div}_{L} (u_1 ,..., u_{2n})=\sum_{j=1}^n (X_j u_j +Y_j
u_{n+j}),$$ $\{X_j, Y_j\}_{j=1, ... ,n}$ are the Greiner type vector
fields (1.1) for general $\ k\geq 1.$
 They
obtained a fundamental solution for $L_{p,k}$ at the origin for
$1<p<\infty,$
$$
\Gamma_p
=\begin{cases} C_{p, k} \ \rho^{\frac{p-Q}{p-1}}\ , & \quad p\neq Q\\
 C_{Q, k} \log \frac{1}{\rho}\ , & \quad p=Q
\end{cases};
$$
where $\rho(z,t)=(|z|^{4k}+t^2)^{1/4k}, Q=2n+2k.$

\begin{remark} Note that when $p = 2$ and $ k = 1,\ L_{p,k}$ becomes the
sub-Laplacian $\Delta_{\mb H^n}$ on the Heisenberg group $\mb H^n.$
If $p=2$ and $k =2, 3, ...$, $L_{p,k}$ is a Greiner  operator (see
 \cite{BGG-JMPA-1998}, \cite{Greiner-cjm79}). Also we note that vector fields in (1.1) do
not possess the translation invariance and they do not satisfy
H\"{o}rmander's condition for $ k>1, k\notin Z.$ Finally we mention
that $L_{p,k}u =0$ is the Euler-Lagrange equation associated to the
functional
$$\int|\na_{L} u|^p, \ \ \  p>1$$
for  functions $u$ satisfying $u, \na_{L} u \in L^p.$
\end{remark}

In the present paper we introduce a family of the vector fields
$X=\{X_1 , ..., X_m\}$ and the corresponding $p-$sub-Laplacian on
H-type groups  generalizing  both of the works above. We let \begeqn
X_j =\partial_j +\frac{1}{2}k|z|^{2k-2} \partial_{[z, e_j]}, \
j=1,2...,m,\endeqn where $\partial_{j}=\partial_{e_j}$,
$\partial_{[z, e_j]}$ are the directional  derivatives,
$\{e_j\}_{j=1,...,m}$ is an orthonormal basis of  $ V$ and $ k\geq 1$
is a fixed parameter.
 We consider the corresponding degenerate p-Laplacian operator
 \begeqn L_{p, k} u= \text{div}_X
(|\na_{X} u|^{p-2} \na_X u),\endeqn
 where
  $$ \na_{X}u =(X_1 u ,..., X_m u), \ \ \text{div}_X
(u_1 ,..., u_{m})=\sum_{j=1}^m X_j u_j. $$

 A natural family of anisotropic dilations
attached to $L_{p, k}$ is
\begeqn\delta_{\lambda}: (z, t) \mapsto (w, s):= (\lambda
z, \lambda^{2k} t),\  \lambda>0,\ (z, t)\in \mb G=\mb
R^{m+q}.\endeqn
It is easy to verify that volume
is transformed by $\delta$ via
$$dw ds= \lambda^Q dz
dt,
$$ where
 $$
Q := m + 2kq,$$ which we may call the degree of homogeneity and  is
the homogeneous dimension in the case $k=1$. We define a
corresponding {\it homogeneous norm} by \begeqn \label{def-d} d(z,
t) := (|z|^{4k} + 16 |t|^{2})^{1/4k}.\endeqn

\section{Fundamental Solutions}
 The main result of this section is the following

\begin{theorem}
Let $\mb G$ be a  H-type group
identified with its Lie algebra $\mathfrak G$
as in (\ref{Lie-prod}).
Then for
$1<p<\infty$,
$$
\Gamma_p
=\begin{cases} C_p \ d^{\frac{p-Q}{p-1}}\ , & \quad p\neq Q\\
 C_Q \log \frac{1}{d}\ , & \quad p=Q
\end{cases}
$$ is a fundamental
solution of $L_{p,k}$  with singularity at the identity element
$0\in \mb G$.
Here $d(z,t)$ is defined in (\ref{def-d}),
 $$ C_p
=\frac{p-1}{p-Q}(\sigma_p)^{-\frac{1}{p-1}},\ \ \
C_Q=-(\sigma_Q)^{-\frac{1}{Q-1}},$$ and
$$\sigma_p =\l(\frac{1}{4}\r)^{q-\frac{1}{2}}\frac{\pi^{\frac{q+m}{2}}\Gamma(\frac{(2k-1)p+m}{4k})}
{\Gamma(\frac{m}{2})\Gamma(\frac{(2k-1)p+Q}{4k})}.$$
\end{theorem}

We prove first some technical identities, which might be of
independent interests.
\begin{lemma} Let $\epsilon >0$ and
$d_{\varepsilon}=(d^{4k}+\gep^{4k})^{\frac{1}{4k}}$. Then we have
 \begeqn
|\nabla_ X d_{\gep}|^2=
\sum_{j=1}^m |X_j (d_{\gep})|^2  = \frac{d^{4k}}{d_{\gep}^{8k-2}}|z|^{4k-2},
 \endeqn
\begeqn
L_{2, k}d_{\gep}^{4k}
=\sum_{j=1}^m X_j^2 (d_{\gep}^{4k})
 = 4k(4k-2+Q)|z|^{4k-2},
\endeqn
and
 \begeqn
L_{2, k} d_{\gep} = \sum_{j=1}^m X_j^2 d_{\gep} =|\na_X d_\gep|^2
\frac{d_{\gep}^{4k-1}}{d^{4k}} \l\{4k+Q-2
-(4k-1)d_{\gep}^{-4k}d^{4k}\r\}.
\endeqn
\end{lemma}
\begin{proof}
By direct computation,
 \begeqn
 \bs
X_j (d_{\gep})& = \frac{1}{4k}d_{\gep}^{1-4k} X_j
(d_{\gep}^{4k})\\
&= \frac{1}{4k}d_{\gep}^{1-4k} \left[4k|z|^{4k-2}  \langle z, e_j
\rangle + 16 k|z|^{2k-2}  \langle t, [z, e_j] \rangle
\right]\\
&= d_{\gep}^{1-4k}\left[|z|^{4k-2} \langle z, e_j \rangle +
4|z|^{2k-2} \langle J_t (z), e_j \rangle \right].
 \end{split}
 \endeqn
 However
  \begeqn\langle J_t (z), z\rangle=\langle t , [z,z] \rangle =0, \ \ \
 \langle J_t (z), J_t (z)\rangle =|t|^2 |z|^2,\endeqn
thus
 $$\sum_{j=1}^m \langle z, e_j \rangle \langle J_t (z), e_j \rangle =\langle J_t (z), z \rangle=0.$$
Consequently
 \begeqn
 \bs
|\na_X d_\epsilon|^2 &=
\sum_{j=1}^m |X_j (d_{\gep})|^2  =d_{\gep}^{2-8k}\left[|z|^{8k-4} |z|^2  +
16|z|^{4k-4} |t|^2 |z|^2 \right]\\
&= \frac{d^{4k}}{d_{\gep}^{8k-2}}|z|^{4k-2},
 \end{split}
 \endeqn
proving the first identity.
Continuing the previous computation of $X_j d_\epsilon$, we find
 \begeqn\bs \sum_{j=1}^m X_j^2 (d_{\gep}^{4k})
&=\sum_{j=1}^m X_j
\l[X_j (d^{4k})\r]\\
&=\sum_{j=1}^m X_j \l[4k \l(|z|^{4k-2} \langle z, e_j \rangle +
4|z|^{2k-2}
\langle J_t (z), e_j \rangle \r)\r]\\
&=4k \sum_{j=1}^m \l\{(2k-1)|z|^{4k-4} 2 \langle z, e_j \rangle^2 +
|z|^{4k-2} \r. \\
& \ \ \ \ \ \ \l. + 8(k-1)|z|^{2k-4}\langle z, e_j \rangle \langle
J_t (z), e_j \rangle + 2k |z|^{4k-4}\langle J_{[z, e_j]} (z), e_j
\rangle \r\}.
\end{split}
\endeqn
To compute the last term in (3.13), we choose an orthonormal basis
$\{t_i\}_{i=1,...q}$ of $\mf t,$ then
 \begeqn\bs \sum_{j=1}^m \langle J_{[z, e_j]}(z),
e_j \rangle &=\sum_{j=1}^m |[z, e_j]|^2 = \sum_{j=1}^m \sum_{i=1}^q
\langle t_i , [z, e_j] \rangle^2 =\sum_{i=1}^q \sum_{j=1}^m \langle
J_{t_i}(z) , e_j \rangle^2\\
&=\sum_{i=1}^q |t_i|^2 |z|^2 =q |z|^2.
\end{split}
\endeqn
Therefore
\begeqn\bs \sum_{j=1}^m X_j^2 (d_{\gep}^{4k})
 &= 4k \l\{(4k-2)|z|^{4k-2}+ m |z|^{4k-2}
+2k|z|^{4k-4}\cdot q |z|^2
\r\}\\
& = 4k(4k-2+Q)|z|^{4k-2},
\end{split}
\endeqn
where $Q=m+2kq.$
We can find $X_j^2 d_\epsilon$
in terms of $X_j^2 d_\epsilon^{4k}$
and $|X_j^2 d_\epsilon|^2$. Indeed
 \begeqn X_j^2 (d_{\gep}^{4k})=X_j (4k
d_{\gep}^{4k-1} X_j d_{\gep})=4k d_{\gep}^{4k-1} X_j^2 d_{\gep } +
4k(k-1)d_{\gep}^{4k-2}|X_j d_{\gep}|^2,\endeqn
thus  \begeqn\bs \sum_{j=1}^m X_j^2 d_{\gep}
&=\frac{1}{4k}d_{\gep}^{1-4k}\l\{\sum_{j=1}^m X_j^2
(d_{\gep}^{4k})-4k(k-1)d_{\gep}^{4k-2}\sum_{j=1}^m |X_j d_{\gep}|^2\r\}\\
&=\frac{1}{4k}d_{\gep}^{1-4k}\l\{4k(4k+Q-2)|z|^{4k-2}
-4k(4k-1)d_{\gep}^{-4k}d^{4k}|z|^{4k-2}\r\}\\
&=d_{\gep}^{1-4k} |z|^{4k-2}\l\{4k+Q-2
-(4k-1)d_{\gep}^{-4k}d^{4k}\r\}\\
&= |\na_X d_\gep|^2 \frac{d_{\gep}^{4k-1}}{d^{4k}} \l\{4k+Q-2
-(4k-1)d_{\gep}^{-4k}d^{4k}\r\},
\end{split}\endeqn
by using the first identity.
\end{proof}

We prove now Theorem 1.
\begin{proof}
 We consider the case $1<p<Q$ first.  Denote
$d_{\varepsilon}=(d^{4k}+\gep^{4k})^{\frac{1}{4k}},\ \
\varepsilon>0.$ We compute $L_{p,k}
(d_{\varepsilon}^{\frac{p-Q}{p-1}}).$ The function
$v=d_{\varepsilon}^{\frac{p-Q}{p-1}}$ is of the form $v=f\circ
d_{\varepsilon}$ with $f(x)=x^{\frac{p-Q}{p-1}}.$ For $f\in C^2 (\mb
R^{+}),$ we have \begeqn\bs L_{p,k} (f\circ d_{\varepsilon}) &=f'|f'|^{p-2}|\na_{X} d_{\gep}|^{p-2} \sum_{j=1}^m X_j^2 d_{\gep}+
|\na_{X} d_{\gep}|^{p-2} \sum_{j=1}^m X_j d_{\gep}\cdot X_j \l(f'|f'|^{p-2}\r)\\
&\ \ \ +f'|f'|^{p-2} \sum_{j=1}^m X_j d_{\gep}\cdot X_j \l(|\na_{X} d_{\gep}|^{p-2}\r)\\
&= I_1 +I_2 +I_3.
\end{split}\endeqn
$I_1$ and $I_2$ can be found by using the Lemma 1,
 \begeqn\bs I_1 &= f'|f'|^{p-2}|\na_{X} d_{\gep}|^{p-2} |\na_{X} d_{\gep}|^{2}
\frac{d_{\gep}^{4k-1}}{d^{4k}}\l\{4k+Q-2
-(4k-1)d_{\gep}^{-4k}d^{4k}\r\}\\
&= f'|f'|^{p-2}|\na_{X} d_{\gep}|^{p}
\l\{(4k+Q-2)\frac{d_{\gep}^{4k-1}}{d^{4k}}
-\frac{(4k-1)}{d_{\gep}}\r\},
\end{split}\endeqn
\begeqn\bs I_2 &= |\na_{X} d_{\gep}|^{p-2}\sum_{j=1}^m X_j
d_{\gep}\cdot \l\{ f''|f'|^{p-2}X_j d_{\gep} + (p-2)|f'|^{p-2}f''
X_j
d_{\gep}\r\}\\
&=|\na_{X} d_{\gep}|^{p} \l\{ f''|f'|^{p-2} + (p-2)|f'|^{p-2}f''\r\}\\
&=|f'|^{p-2}|\na_{X} d_{\gep}|^{p}\l\{ (p-1)f''\r\}.
\end{split}\endeqn
Using $X_j |\nabla_X d_\epsilon |^{p-2}= \frac{p-2}2 |\nabla_X
d_\epsilon|^{p-4}X_j |\nabla_X d_\epsilon |^2$ and the Lemma 1, we
find
  \begeqn\bs I_3 &= f'|f'|^{p-2}\sum_{j=1}^m X_j d_{\gep}\cdot
\frac{p-2}{2}|\na_{X} d_{\gep}|^{p-4} X_{j}(|\na_{X} d_{\gep}|^{2})\\
&=\frac{p-2}{2}f'|f'|^{p-2}|\na_{X} d_{\gep}|^{p-4}\sum_{j=1}^m X_j
d_{\gep}\cdot X_j \l(d_{\gep}^{2-8k}d^{4k}|z|^{4k-2}\r)\\
&=\frac{p-2}{2}f'|f'|^{p-2}|\na_{X} d_{\gep}|^{p-4}\sum_{j=1}^m X_j
d_{\gep}\cdot\l\{ (2-8k)d_{\gep}^{1-8k}d^{4k}|z|^{4k-2}
X_{j}d_{\gep} \r. \\
&\ \ \ \ \ \ \ \  \ \ \ \  \ \ \ \ \ \ \ \ \ \ \ \ \ \  \ \ \ \ \ \
\ \ \ \ \  \ + 4k
d_{\gep}^{2-8k}d^{4k-1}|z|^{4k-2} X_{j}d \\
&\ \ \ \ \ \ \ \  \ \ \ \  \ \ \ \ \ \ \ \ \ \ \ \  \ \ \ \ \ \ \ \
\ \ \ \ \ \  \l.+ (4k-2)
d_{\gep}^{2-8k}d^{4k}|z|^{4k-4} \langle z, e_j \rangle \r\}\\
&=\frac{p-2}{2}f'|f'|^{p-2}|\na_{X}
d_{\gep}|^{p-4}\l\{(2-8k)d_{\gep}^{3-16k}d^{8k}|z|^{8k-4}+ 4k
d_{\gep}^{3-12k}d^{4k}|z|^{8k-4} \r.\\
& \ \ \ \ \ \ \ \ \ \ \ \ \ \ \ \ \ \ \ \ \ \ \ \ \ \ \ \ \ \ \ \ \
\ \ \ \ \ \  \l.+(4k-2)
d_{\gep}^{3-12k}d^{4k}|z|^{8k-4} \r\}\\
&=(p-2)(4k-1)f'|f'|^{p-2}|\na_{X}
d_{\gep}|^{p}\frac{\gep^{4k}}{d_{\gep}\ d^{4k}}.
\end{split}\endeqn
Hence,
 \begin{equation}\bs   &\ \ \ \ \  L_{p,k} (f\circ d_{\varepsilon})\\&= I_1 +I_2 +I_3\\
&=|f'|^{p-2}|\na_{X}
d_{\gep}|^{p}\l\{(p-1)f''+f'\l[\frac{(Q-1)d^{4k}+(4kp-4k+Q-p)\gep^{4k}}{d_{\gep}
d^{4k}}\r]\r\}.
\end{split}\end{equation}
Taking $f(x)=x^{\frac{p-Q}{p-1}} \ \ (x>0)$ the above is
\begeqn\bs   L_{p,k} \l(d_{\gep}^{\frac{p-Q}{p-1}}\r)
&=\l|\frac{p-Q}{p-1}d_{\gep}^{\frac{1-Q}{p-1}}\r|^{p-2}
\l(\frac{d^{2k}|z|^{2k-1}}{d_{\gep}^{4k-1}}\r)^{p}\l\{\frac{p-Q}{p-1}(1-Q)d_{\gep}^{\frac{2-p-Q}{p-1}}\r.\\
&\ \ \ \ \ \ \ \ \ \  \ \ \ \ \
\l.+\frac{p-Q}{p-1}d_{\gep}^{\frac{1-Q}{p-1}}\l[\frac{(Q-1)d^{4k}+(4kp-4k+Q-p)\gep^{4k}}{d_{\gep}
d^{4k}}\r]\r\}\\
&=-\l(\frac{Q-p}{p-1}\r)^{p-1}d_{\gep}^{1-Q}\l(\frac{d^{2k}|z|^{2k-1}}{d_{\gep}^{4k-1}}\r)^{p}
\l\{(4kp-4k+Q-p)\frac{\gep^{4k}}{d^{4k}d_{\gep}}\r\}\\
&=-\l(\frac{Q-p}{p-1}\r)^{p-1}(4kp-4k+Q-p)\frac{d^{2kp-4k}|z|^{(2k-1)p}\gep^{4k}}{d_{\gep}^{(4k-1)p+Q}}\\
&=\gep^{-Q}\psi(\delta_{1/\gep}(z,t)),
\end{split}\endeqn
where
$$
\psi(z,t):=-\l(\frac{Q-p}{p-1}\r)^{p-1}(4kp-4k+Q-p)\frac{d^{2kp-4k}|z|^{(2k-1)p}}{(1+d^{4k})^{(4kp-p+Q)/4k}}.$$

Now for any $\varphi\in C_0^{\infty}(\mb G),$ it follows that \begeqn\bs
\langle L_{p,k} (d^{\frac{p-Q}{p-1}}), \varphi \rangle
&=\lim_{\gep\rightarrow 0}\int_{\mathbb{G}}L_{p,k}
(d_{\gep}^{\frac{p-Q}{p-1}})\varphi\\
&=\lim_{\gep \rightarrow
0}\gep^{-Q}\int_{\mb G}\psi(\delta_{1/\gep}(z,t))\varphi(z,t)\\
&=\lim_{\gep \rightarrow 0}\int_{\mb G}\psi(z,t)\varphi(\gep z, \gep^{2k}t)\\
&=\varphi(0)\int_{\mb G}\psi(z,t).
\end{split}\endeqn

Finally we evaluate the integral $\int_{\mathbb{G}}\psi(z,t).$ We
use the polar coordinates $z=r z^\ast$ with  $r=d$
and
$z^\ast \in S:=\{g\in
\mb G: \ d(g)=1 \}$ being
the sphere with respect to $d$ 
. By a general integral formula
 on homogeneous groups (see
\cite{Folland-stein1982Hardyspace}, Proposition 1.15)
we have
\begin{equation*}
\begin{split}
 &
-\int_{\mathbb{G}}\psi(z,t)\\
&=
(4kp-4k+Q-p)\int_{\mb G}\frac{d^{2kp-4k}|z|^{(2k-1)p}}{(1+d^{4k})^{(4kp-p+Q)/4k}}\\
&=(4kp-4k+Q-p)\int_{S}|z^*|^{(2k-1)p}\int_{0}^{\infty}\frac{r^{-4k-1}}{(1+r^{-4k})^{(4kp-p+Q)/4k}}drd\sigma
\\&=(4kp-4k+Q-p)\int_{S}|z^*|^{(2k-1)p}d\sigma
\frac{1}{4k}\int_{1}^{\infty}t^{\frac{p-Q-4kp}{4k}}dt\\
&=\int_{S}|z^*|^{(2k-1)p}d\sigma,
\end{split}\end{equation*}

 Denote temporarily
$\gga =(2k-1)p$. We
use the usual trick to evaluate the integral on the sphere,
 replacing it by an integral on the ball,
 \begin{equation*}\bs \int_{S}|z^*|^{\gga}d\sigma
&=(Q+\gga)\int_0^1 r^{\gga
+Q-1}dr\int_{S}|z^*|^{\gga}d\sigma\\
&=(Q+\gga)\int_{S}\int_0^1 |rz^*|^{\gga} r^{Q-1}drd\sigma\\
&=(Q+\gga)\int_{d<1}|z|^\gga,
\end{split}
\end{equation*}
and furthermore
 \begin{equation*}\bs \int_{d<1}|z|^\gga &=\int_{|t|<\frac{1}{4}}\int_{|z|<(1-16|t|^2)^{\frac{1}{4k}}}|z|^{\gga}dzdt\\
&=\omega_{m-1}\int_{|t|<\frac{1}{4}}\int_{0}^{(1-16|t|^2)^{\frac{1}{4k}}}r^{\gga+m-1
}drdt\\
&=\frac{\omega_{m-1}\omega_{q-1}}{\gga
+m}\int_0^{\frac{1}{4}}(1-16s^2)^{\frac{\gga +m}{4k}}s^{q-1}ds\\
&=\frac{\omega_{m-1}\omega_{q-1}}{2(\gga +m)}\l(\frac{1}{4}\r)^q
\int_0^1 (1-\rho)^{\frac{\gga +m}{4k}}\rho^{\frac{q-2}{2}}d\rho\\
&=\frac{\omega_{m-1}\omega_{q-1}}{2(\gga +m)}\l(\frac{1}{4}\r)^q
\frac{\Gamma(\frac{\gga
+m+4k}{4k})\cdot\Gamma(\frac{q}{2})}{\Gamma(\frac{\gga +m+4k+2kq}{4k})}\\
&=\frac{1}{2(\gga+Q)}\l(\frac{1}{4}\r)^{q-1}\frac{\pi^{\frac{q+m}{2}}\cdot
\Gamma(\frac{\gga+m}{4k})}{\Gamma(\frac{m}{2})\cdot\Gamma(\frac{\gga+Q}{4k})}.
\end{split}\end{equation*}
Thus, \begin{equation*}\int_{S}|z^*|^{(2k-1)p}d\sigma
=\l(\frac{1}{4}\r)^{q-\frac{1}{2}}\frac{\pi^{\frac{q+m}{2}}\cdot
\Gamma(\frac{(2k-1)p+m}{4k})}{\Gamma(\frac{m}{2})\cdot\Gamma(\frac{(2k-1)p+Q}{4k})}, 
\end{equation*}
and substituthing this into the previous formula
for 
$
-\int_{\mathbb{G}}\psi(z,t)
$
we find
$$\int_{\mb G}\psi(z,t)=-\l(\frac{Q-p}{p-1}\r)^{p-1}
\l(\frac{1}{4}\r)^{q-\frac{1}{2}}\frac{\pi^{\frac{q+m}{2}}\cdot
\Gamma(\frac{(2k-1)p+m}{4k})}{\Gamma(\frac{m}{2})\cdot\Gamma(\frac{(2k-1)p+Q}{4k})}$$
proving Theorem 1  for $1<p<Q.$

An direct examination shows that the formula also holds for $p>Q,$
and the critical case $p=Q$ can be treated similarly, we omit the
details.
\end{proof}

By a similar method as in Theorem 1, we can also obtain a
fundamental solution for a class of  weighted p-Laplacian operators
on the H-type group $\mb G =\mb R^m \oplus \mb R^q,$ \begeqn
L_{p,k,w}= \text{div}_X (|\na_{X} u|^{p-2}w  \na_{X} u),\endeqn
$$\l(w=d^{\ga}|\na_{X} d|^{\gb},\ \ga
>-m-2kq,  \gb>\max\l\{\frac{1-Q}{4k-1}, -\frac{m}{2k-1}-1\r\}\r)$$
 where $\{X_j\}_{j=1, ..., m }$ is taken from (2.3)
and $d(z, t)$  from (2.6).
\begin{theorem}Let $\mb G$ be the H-type group above
 $L_{p,k,w}$ the $p$-sub-Laplacian defined defined as in (3.25). Then for
$1<p<\infty$
$$ \Gamma_{p,w}
=\begin{cases} C_{p,w} \ d^{\frac{p-Q-\ga}{p-1}}\ , & \quad p\neq Q+\ga\\
 C_{Q+\ga , w} \log \frac{1}{d}\ , & \quad p=Q+\ga
\end{cases};
$$
 is a fundamental
solution of $L_{p,k,w}$ with singularity at the identity element
$0\in \mb G$, where
$$  C_{p, w}
=\frac{p-1}{p-Q-\ga}(\sigma_{p, \gb})^{-\frac{1}{p-1}}, \ \
C_{Q+\ga, w}=-(\sigma_{Q+\ga, \gb})^{-\frac{1}{Q+\ga-1}},$$ and
$$\sigma_{p,\gb}=\l(\frac{1}{4}\r)^{q-\frac{1}{2}}\frac{\pi^{\frac{q+m}{2}}}{\Gamma(\frac{m}{2})}
\frac{\Gamma(\frac{(2k-1)(p+\gb)+m}{4k})}{\Gamma(\frac{(2k-1)(p+\gb)+Q}{4k})}.$$
\end{theorem}

\section{Hardy type inequality }

We recall that the classical Hardy inequality
states that, for $n\geq3,$
\begin{equation} \int_{\mb R^n} |\nabla \gp(x)|^2 dx\geq
\left(\frac{n-2}{2}\right)^2 \int_{\mb R ^n}\frac{|
\gp(x)|^2}{|x|^2}dx,
\endeqn
where $\gp \in C_0 ^\infty(\mb R^n\setminus \{0\}).$ 
It
can also be rewritten in terms of certain Schr\"o{}dinger operator.
The inequality and their generalizations are
thus of interests in
 the study of spectral theory of linear
and nonlinear partial differential equations
(see
e. ~g.~ \cite{Peral-JDE1998}, \cite{Goldstein-Kombe-IJEE2005},
\cite{Goldstein-Zhang.Qi.S-JFA2001}).

 In \cite{Garofalo-Lanconelli-Fourier1990} Garofalo and Lanconelli
established the following Hardy inequality on the Heisenberg  group
${\mb H}=\mb H^n$ associated with left-invariant horizontal gradient
$\na_\mb H,$ \begeqn \int_{{\mb H}} |\nabla_{{\mb H}} \gp|^2
dzdt\geq\left(\frac{Q-2}{2}\right)^2\int_{{\mb H}}
\left(\frac{|z|^2}{|z|^4+t^2}\right) |\gp|^2 dzdt,
\endeqn
where $\gp \in C_0 ^\infty({\mb H}\setminus \{0\}),$ $Q=2n+2$ is the
homogeneous dimension of ${\mb H},$ and $\na_{{\mb H}}\gp =(X_1\gp,
X_2\gp,\ccd, X_n\gp, Y_1\gp, \ccd, Y_n\gp), \ X_j=\frac{\pa}{\pa
x_j}+2y_j\frac{\pa}{\pa t},\ Y_j=\frac{\pa}{\pa
y_j}-2x_j\frac{\pa}{\pa t},$ for $ (z,t)\in \mb H, z=(x,y)\in \mb
R^n\times\mb R^n,\ t\in \mb R.$ The $L^p$ version of the inequality
(4.27) has been obtained,  among others,  by Niu, Zhang, Wang in
\cite{NZWpams2001Hardyinequality}, which states that for $1<p<Q$:
\begeqn \int_\mb H |\na_{\mb H} \gp|^p \geq \l(\frac{Q-p}{p}\r)^p
\int_\mb H \l(\frac{|z|}{d}\r)^{p}\frac{|\gp|^p}{d^p}.
\endeqn

In this section we obtain  a Hardy type inequality associated
with the non-invariant vector fields $X=\{X_j\}$ in (2.3) on the
H-type groups by applying the result in Section 3. The inequality in
the present paper might be useful in  eigenvalue problems and
Liouville type theorems for weighted p-Laplacian equation, which we
plan to pursue in some subsequent work. Recall the norm $d$  in
(2.6).
 \begth
Let $\mb G$ be the  H-type group with the homogeneous dimension
$Q=m+2kq$ and   $\ga \in \mb R,
1<p<Q+\ga$
Then the
following inequality holds
for
$ \gp \in C_0^\infty (\mb G \backslash \{0\})$,
\begeqn \int_\mb G d^\ga |\na_{X} \gp|^p
\geq \l(\frac{Q+\ga-p}{p}\r)^p \int_\mb G d^\ga
\l(\frac{|z|}{d}\r)^{(2k-1)p}\l|\frac{\gp}{d}\r|^p.
\endeqn
 Moreover, the
constant $(\frac{Q+\ga-p}{p})^p$ is sharp.
 \endth

In view of the first equality in Lemma 1 (for $\epsilon=0$),
namely $|\nabla_X d|=(\frac{|z|}{d})^{2k-1}$, the above
inequality can also be written as
\begin{equation*}
 \int_\mb G d^\ga |\na_{X} \gp|^p
\geq 
\l(\frac{Q+\ga-p}{p}\r)^p 
\int_\mb G 
d^{\ga -p}
|\nabla_X d|^p |{\gp}|^p.
\end{equation*}

\begin{remark}
If  $q=1$ and $ \ga =0,$ then our Theorem 3 is actually the Theorem
3.1 in \cite{ZN2003Nonlinear-Analysis}.
\end{remark}
For the proof of Theorem 3, we need the following Lemma; see also
\cite{NZWpams2001Hardyinequality} for the case $w=1.$
 \beglemm
Let $w\geq 0$ be a weight function in $\gO \subset \mb G$ and
$L_{p,k,w}u=\text{div}_X (|\na_{X} u|^{p-2}w \na_{X} u).$  Suppose
that for some $\gl
>0,$ there exists $ \ v\in C^\infty(\gO), v>0$ such that
\begeqn -L_{p,k,w}v \geq \gl g v^{p-1}
\endeqn
for some $g\ge 0$,  in the sense of distribution
acting on non-negative test functions.
 Then for any
$u\in HW_0^{1,p}(\gO,w),$ it holds that
$$\int_{\gO}|\na_{X}u|^p w \geq \gl \int_\gO g|u|^p,$$
where  $HW_0^{1,p}(\gO, w)$ denote the closure of
$C_0^{\infty}(\gO)$ in the norm $(\int_{\gO}|\na_{X}u|^p
w)^{\frac{1}{p}}.$
\endlemm

\begpr We take $\frac{\gv^p}{v^{p-1}}$ as a test function in (4.30),
where $\gv\in C_0^\infty(\gO), \gv\geq 0$,
$$
I:=\int_\gO w |\na_{X} v|^{p-2}\na_{X} v \cdot\na_{X} \l(\frac{\gv^p}{v^{p-1}}\r)
 \geq \gl \int_\gO g \gv^p.$$
We shall prove \begeqn \int_{\gO}w|\na_{X}\gv|^p - I \ge 0 \endeqn
 which together with the previous inequality implies Lemma 2 for $u=\gv \in C_0^{\infty}(\Omega)$.
Now, the above is an integration with integrand (disregarding the
common factor $w$),
 \begeqn
 \bs
&\ \ \ \ \  |\na_{X}\gv|^p- |\na_{X} v|^{p-2}\na_{X}
\l(\frac{\gv^p}{v^{p-1}}\r)\cdot
 \na_{X} v
\\
 &= |\na_{X}\gv|^p-p\frac{\gv^{p-1}}{v^{p-1}}|\na_{X} v|^{p-2}
 \na_{X}\gv\cdot\na_{X}v
  +(p-1)\frac{\gv^{p}}{v^{p}}|\na_{X}
 v|^{p}
\\
 & = \frac{1}{v^p}\l(v^p|\na_{X}\gv|^p+(p-1)\gv^p|\na_{X}v|^p
 -pv\gv^{p-1}|\na_{X}v|^{p-2}\na_{X}\gv\cdot\na_{X}v\r).
 \end{split}
 \endeqn
We estimate last term from above
using the  Young's inequality
$$ ab\leq
\frac{1}{p}\ a^p+(1-\frac{1}{p})\ b^{\frac{p}{p-1}},$$
and get
 \begeqn
 \bs pv\gv^{p-1}|\na_{X}v|^{p-2}\na_{X}\gv\cdot\na_{X}v &\leq p v|\na_{X}\gv|\cdot
\gv^{p-1}|\na_{X}v|^{p-1}\\&\leq p\ \left[\frac{v^p|\na_{X}
\gv|^p}{p}+\frac{p-1}{p}\gv^p|\na_{X} v|^p\right]\\
&= v^p|\na_{X} \gv|^p+(p-1)\gv^p|\na_{X} v|^p.
\end{split}
 \endeqn
  Hence (4.31) follows. The proof of Lemma 2 is finished by taking
 $\gv\rightarrow u.$
\endpr
We prove now Theorem 3.
\begin{proof}\ \  Case (i): $p\neq Q.$
We claim that  the conditions in Lemma 2 are satisfied with
$$
w=d^{\ga}, \quad v=d^{\frac{p-Q-\ga}{p}},
\quad  g=d^\ga\frac{|z|^{(2k-1)p}}{d^{2kp}},\quad
\lambda=\l(\frac{Q+\ga-p}{p}\r)^{p},
\quad \Omega=\mb G\backslash\{0\},
$$
which then proves the Theorem. Indeed, for  any $\gv\in
C_0^\infty (\mb G\backslash \{0\})$ we have
\begeqn \bs   \langle
-L_{p,k,w}v, \gv\rangle
&=-\left(\frac{Q+\ga-p}{p}\right)^{p-1}\int_{\mb G}
(d^{\frac{Q+\ga}{p}-Q}|\na_{X} d|^{p-2} \na_{X} d)\cdot
\na_{X}\gv\\
&=-\left(\frac{Q+\ga-p}{p}\right)^{p-1}\int_{\mb G} (d^{1-Q}|\na_{X}
d|^{p-2} \na_{X} d)\cdot d^{\frac{Q+\ga-p}{p}}
\na_{X}\gv\\
&=-\left(\frac{Q+\ga-p}{p}\right)^{p-1}\int_{\mb G} (d^{1-Q}|\na_{X}
d|^{p-2} \na_{X} d)\cdot \na_{X}(\gv\cdot
d^{\frac{Q+\ga-p}{p}})\\
&\ \ \ \ +\left(\frac{Q+\ga-p}{p}\right)^{p-1}\int_{\mb G}
(d^{1-Q}|\na_{X} d|^{p-2} \na_{X} d)\cdot \na_{X}(d^{\frac{Q+\ga-p}{p}})\ \gv.\\
\end{split}
\end{equation}
Denoting $C_{p,Q}=\left|\frac{p-1}{p-Q}\right|^{p-2}\frac{p-1}{p-Q}$
and rewriting
$$d^{1-Q}|\na_{X} d|^{p-2} \na_{X} d
=C_{p,Q}\l|\na_{X} \l(d^{\frac{p-Q}{p-1}}\r)\r|^{p-2} \na_{X}
\l(d^{\frac{p-Q}{p-1}}\r)$$
we see that
 (4.34) is
 \begeqn \bs  &\ \ \ \  \langle -L_{p,k,w}v, \gv\rangle\\
&=-C_{p,Q}\left(\frac{Q+\ga-p}{p}\right)^{p-1}\int_{\mb G}
\l|\na_{X} \l(d^{\frac{p-Q}{p-1}}\r)\r|^{p-2} \na_{X}
\l(d^{\frac{p-Q}{p-1}}\r)\cdot
\na_{X}\l(\gv d^{\frac{Q+\ga-p}{p}}\r)\\
&\ \ \ \ + \left(\frac{Q+\ga-p}{p}\right)^{p-1}\int_{\mb G}
\l(d^{1-Q}|\na_{X} d|^{p-2} \na_{X} d\r)\cdot
\na_{X}\l(d^{\frac{Q+\ga-p}{p}}\r)
\gv.\\
\end{split}
\endeqn
However
 the first integral in (4.35) is zero by Theorem 2, 
since $\phi$ is supported away from $0$, and 
we find  \begeqn \bs
\langle -L_{p,k,w}v, \gv\rangle
&=\left(\frac{Q+\ga-p}{p}\right)^{p-1}\int_{\mb G} d^{1-Q}|\na_{X}
d|^{p-2} \na_{X} d \cdot \na_{X}(d^{\frac{Q+\ga-p}{p}})
\gv\\
&=\l(\frac{Q+\ga-p}{p}\r)^{p}\int_{\mb G}
d^{\frac{Q+\ga}{p}-1-Q}|\na_{X} d|^{p}\gv\\
&= \l(\frac{Q+\ga-p}{p}\r)^{p}\int_{\mb G}
 d^\ga d^{\frac{p-Q-\ga}{p}(p-1)}\frac{|z|^{(2k-1)p}}{d^{2kp}}\gv\\
&= \l(\frac{Q+\ga-p}{p}\r)^{p}\int_{\mb G}
 d^\ga\frac{|z|^{(2k-1)p}}{d^{2kp}} v^{p-1}\gv,
\end{split}
\endeqn
where in the second last equality we have used
Lemma 1 that $|\na_{X} d|^p=\l(\frac{|z|}{d}\r)^{(2k-1)p}$. This proves our claim.

 Case (ii): $p=Q.$

 The proof is almost the
same as the above  once we notice the following  fact: $C_{Q}\log
\frac{1}{d}$ is a fundamental solution of $L_{Q,k}$ on $\mb G,$ and
$$d^{1-Q}|\na_{X} d|^{Q-2} \na_{X} d =-|\na_{X} \log
(d^{-1})|^{Q-2} \na_{X} \log (d^{-1}).$$

It remains to show  the sharpness of the constant $(\frac{Q+\ga-p}{p})^p$.
This is equivalent to show that  any constant
$B>0$ for which
 the inequality
\begin{equation}
\label{B-const}
 \int_\mb G d^\ga |\na_{X} \gp|^p
\ge
B \int_\mb G 
d^{\ga -p}
|\nabla_X d|^p |{\gp}|^p
\end{equation}
holds must satisfy $B\le (\frac{Q+\ga-p}{p})^p$.
We shall construct a sequence  $\{u_j\}_{j=1}^\infty$
of functions so that the inequality (4.29) approximates
to an identity up to the order $O(1)$ in $j$.
 Given any positive integer $j$
it is elementary that
 there exists $\psi_j$ in $C^\infty_0(0, \infty)$ 
such that
 $\text{supp}\, \phi_j =[2^{-j-1}, 2]$, 
 $\psi_j(x)=1$ on
$[2^{-j}, 1]$, 
and $|\psi_j^\prime (x)|\le C 2^{j}$
on $[2^{-j-1}, 2^{-j}]$, where
$C$ is a constant independent of $j$.
Let
$$
u_j (z, t)=  
 d(z, t)^{\frac{p-Q-\ga}{p}-\frac 1j}\,
\psi_j(d(z, t)).
$$
 Clearly $u_j \in C^\infty(G\setminus \{0\})$ and
is radial.
The gradient is given by
\begin{equation}
\label{nabla-u}
 \na_{X} u_{j}
=
\begin{cases} 0, & 0\leq d<2^{-j-1}, \,\text{or}\, \,
d>2
 \\
- ({\frac{Q+\ga-p}{p}+\frac 1j})
 d^{-\frac{Q+\ga+p\frac 1j}{p}}\ \na_{X}d, & 2^{-j} < d< 1
\end{cases}
\end{equation}
The left hand side of the above inequality is
$$
LHS= \int_\mb G
=\int_{2^{-j} < d <1}  +
\int_{2^{-j-1} < d <2^{-j}} +  
\int_{1 < d <2} 
=\int_{2^{-j} < d <1}  + I + II.
$$
The first integration is
$$
\int_{2^{-j} < d <1} 
d^{\alpha} |\nabla_X u_j|^p
=
 ({\frac{Q+\ga-p}{p}+\frac 1j})^p
\int_{2^{-j} < d <1} 
 d^{-\frac{Q+\ga+p\frac 1j}{}}\
| \na_{X}d|^p.
$$
This  can  be computed by using the
polar coordinates as in proof of Theorem 1
and is
$$ ({\frac{Q+\ga-p}{p}+\frac 1j})^p
C_0 
j,
$$
where $C_0
=\frac {({2^p-1})}p\int_{S} |z|^{p(2k+1)}
$  (and is evaluated
in the proof of Theorem 1). Similarly,
$$
RHS=B\int_{2^{-j} < d <1}  + III + IV.
$$
The first integration is  precisely the same as above 
and is 
$$
B\int_{2^{-j} < d <1}=B\, C_0 j,
$$
with the same constant $C_0$.
It is easy to estimate the error terms
and they are all bounded
$$
I, II, III, IV \le C.
$$
The inequality (\ref{B-const}) now becomes
$$
 ({\frac{Q+\ga-p}{p}+\frac 1j})^p
C_0 j + I +II
\ge BC_0 j + III + IV.
$$
Dividing both sides by $j$ and letting
$j\to \infty$ prove our claim.
\end{proof}

 An immediate consequence of Theorem 3
is the following corollary,  known also as the uncertainty
principle, this can be proved by estimating the left hand side using
H\"{o}lder inequality together with inequality (4.29) for $\ga =0$.
\begin{corollary}Let $\mb G$ be the  H-type group with the homogeneous dimension
$Q=m+2kq$ associated with the dilations (2.5). $ u\in C_0^\infty
(\mb G \backslash\{0\} ), \frac{1}{s}+\frac{1}{t}=1\  (1<s<Q).$ Then
$$\l(\int_{\mb G}|z|^t |u|^t \r)^{\frac{1}{t}}\l(\int_{\mb G}|\na_{X}u|^s \r)^{\frac{1}{s}}
\geq \frac{Q-s}{s}\int_{\mb G}\frac{|z|^{2k}}{d^{2k}}|u|^2.$$
\end{corollary}

\def\cprime{$'$} \def\cprime{$'$}

\bigskip

 Yongyang Jin: Department of Applied Mathematics,
 Zhejiang University of Technology,   Hangzhou, 310032,  China;
yongyang@zjut.edu.cn.

    Genkai Zhang:
Mathematical Sciences, Chalmers University of Technology and Mathematical Sciences, G\"oteborg
University,   G\"oteborg, Sweden;  genkai@math.chalmers.se

\end{document}